\documentclass[12pt]{article}
\usepackage{amsfonts}
\usepackage{amsmath}
\usepackage{amssymb}

\usepackage[T2A]{fontenc} 
\usepackage[russian,english]{babel}
\usepackage{citehack}     

\newcommand\ov{\overline}

\newcommand\opn{\operatorname}

\newcommand\al{\alpha}

\newcommand\ga{\gamma}

\newcommand\De{\Delta}

\renewcommand\phi{\varphi}

\newcommand\ep{\varepsilon}

\def\supl{\sup\limits}

\newlength{\MMM}
\newlength{\MMMM}
\newcommand{\loc}{\opn{loc}}

\newtheorem{theor}{рЕНПЕЛЮ}[section]
\newtheorem{theorem}[theor]{рЕНПЕЛЮ}
\newtheorem{proposition}[theor]{оПЕДКНФЕМХЕ}

\newenvironment{proof}{\par\smallskip{\bf дНЙЮГЮРЕКЭЯРБН.}\rm}
			{\hfill $\square$\medskip}
\newenvironment{proof*nodot}{\par\smallskip{\bf дНЙЮГЮРЕКЭЯРБН}\rm}
			{\hfill $\square$\medskip}

\newenvironment{theorem*}[1]{\par\smallskip{\bf #1.}\it}{\par\medskip}

\begin{document}

\makeatletter


\renewcommand{\subsection}{\@startsection{subsection}{2}%
{0pt}{2.25ex plus 1ex minus .2ex}%
{0.5ex plus .2ex}{\large\bf}
}

\makeatother

\renewcommand{\bibname}{яОХЯНЙ КХРЕПЮРСПШ}

\title{
н ЛСКЭРХОКХЙЮРНПЮУ ХГ ЯНАНКЕБЯЙНЦН ОПНЯРПЮМЯРБЮ $H^\al_p$ Б $H^{-\al}_p$.
}
\author{
л.~х.~мЕИЛЮМ--ГЮДЕ\footnote{пЮАНРЮ ОНДДЕПФЮМЮ ЦПЮМРНЛ
пттх N 01-01-00691}, ю.~ю.~ьЙЮКХЙНБ\footnote{пЮАНРЮ ОНДДЕПФЮМЮ ЦПЮМРНЛ
пттх N 00-15-96100}}


\maketitle

\newcommand{\RR}{\mathbb R}
\newcommand{\unif}{\opn{unif}}
\newcommand{\norm}{|||}
\newcommand{\sD}{\mathcal D}

б ЩРНИ ЯРЮРЭЕ ЛШ АСДЕЛ ХГСВЮРЭ ОПНЯРПЮМЯРБЮ ЛСКЭРХОКХЙЮРНПНБ ХГ
ЯНАНКЕБЯЙНЦН ОПНЯРПЮМЯРБЮ $H_p^\alpha(\RR^n)$ Б ОПНЯРПЮМЯРБН
$H_p^{-\alpha}(\RR^n)$, $1<p<\infty$.
мЮЬЮ ПЮАНРЮ ЪБКЪЕРЯЪ ОПНДНКФЕМХЕЛ ЯРЮРЭХ дФ.~ц.~аЮЙЮ Х ю.~ю.~ьЙЮКХЙНБЮ~\cite{BSh}
Н ЛСКЭРХОКХЙЮРНПЮУ ХГ ОПНЯРПЮМЯРБЮ $H^\al_p$ Б ДСЮКЭМНЕ ОПНЯРПЮМЯРБН $H^{-\al}_p$,
ЦДЕ $p$ Х $p'$ --- ЯНОПЪФЕММШЕ ОН ц\"ЕКЭДЕПС ВХЯКЮ
Р.Е $1/p+1/p'=1$.

лШ ОНКСВХЛ РНВМНЕ НОХЯЮМХЕ ОПНЯРПЮМЯРБЮ ЛСКЭРХОКХЙЮРНПНБ ХГ
$H_p^\alpha(\RR^n)$ Б $H_p^{-\alpha}(\RR^n)$
(ЛШ НАНГМЮВЮЕЛ ЕЦН ВЕПЕГ~$M^\alpha_p$) ОПХ СЯКНБХХ
$\alpha>\min(n/p,n/p')$. б ЯКСВЮЕ $\alpha\le \min(n/p,n/p')$ АСДСР ОНКСВЕМШ
РЕНПЕЛШ БКНФЕМХЪ ЯНАНКЕБЯЙХУ ОПНЯРПЮМЯРБ Я МЕЦЮРХБМШЛХ ХМДЕЙЯЮЛХ
ЦКЮДЙНЯРХ Б ОПНЯРПЮМЯРБЮ~$M^\alpha_p$ (ДНЯРЮРНВМШЕ СЯКНБХЪ
ОПХМЮДКЕФМНЯРХ Й~$M^\alpha_p$).
дНЙЮГЮРЕКЭЯРБЮ МЕЙНРНПШУ ХГ ОПХБНДХЛШУ МХФЕ СРБЕПФДЕМХИ ОНВРХ ДНЯКНБМН ОНБРНПЪЧР
ДНЙЮГЮРЕКЭЯРБЮ ХГ~\cite{BSh}, МН ДКЪ ОНКМНРШ ХГКНФЕМХЪ ЛШ ХУ ОПХБНДХЛ.

бБЕДЕЛ НАНГМЮВЕМХЕ
$$
H_{p,\loc}^\alpha(\RR^n)
=\bigl\{v\in D':
v(x)\phi(x)\in H_p^\alpha \quad \forall\phi\in D\bigr\},
$$
ЦДЕ $D'$~-- ОПНЯРПЮМЯРБН ПЮЯОПЕДЕКЕМХИ МЮ~$D$. тСМЙЖХЧ
$\mu(x)\in H_{p,\loc}^{-\alpha}$ МЮГНБЕЛ
{\it ЛСКЭРХОКХЙЮРНПНЛ\/} ХГ $H_p^\alpha$
Б~$H_p^{-\alpha}$, ЕЯКХ
\begin{equation}\label{II:eq:negmult}
\|\mu(x)\phi(x)\|_{-\alpha,p}
\le C\|\phi\|_{\alpha,p}
\qquad \forall\phi\in D,
\end{equation}
ЦДЕ ВЕПЕГ $C$ НАНГМЮВЮЕРЯЪ (ГДЕЯЭ Х Б ДЮКЭМЕИЬЕЛ)
ЙНМЯРЮМРЮ, МЕ ГЮБХЯЪЫЮЪ НР $\phi\in D$.

б ЩРНЛ НОПЕДЕКЕМХХ ОПЕДОНКЮЦЮЕРЯЪ, ВРН $\phi\in D$ (БЛЕЯРН
$\phi\in H_p^\alpha$), ОНЯЙНКЭЙС ОПХ
$\mu(x)\in H_{p,\loc}^{-\alpha}$ ОПНХГБЕДЕМХЕ
$\mu(x)\phi(x)$ НОПЕДЕКЕМН Б~$D'$ КХЬЭ ДКЪ АЕЯЙНМЕВМН
ЦКЮДЙХУ~$\phi$. нДМЮЙН $D$ ОКНРМН Б $H_p^\alpha$, ОНЩРНЛС
ХГ НЖЕМЙХ~\eqref{II:eq:negmult} ДКЪ $\phi\in D$ ЯКЕДСЕР, ВРН НОЕПЮРНП
СЛМНФЕМХЪ МЮ ТСМЙЖХЧ $\mu(x)$ ЛНФЕР АШРЭ ОПНДНКФЕМ МЮ БЯЕ ОПНЯРПЮМЯРБН
$H_p^\alpha$ ЙЮЙ НЦПЮМХВЕММШИ НОЕПЮРНП ХГ $H_p^\alpha$
Б~$H_p^{-\alpha}$. нВЕБХДМН, ТСМЙЖХХ~$\mu(x)$,
СДНБКЕРБНПЪЧЫХЕ МЕПЮБЕМЯРБС~\eqref{II:eq:negmult},
НАПЮГСЧР АЮМЮУНБН
ОПНЯРПЮМЯРБН, ЙНРНПНЕ ЛШ НАНГМЮВХЛ ВЕПЕГ
$M[H_p^\alpha\to H_p^{-\alpha}]$, ХКХ АНКЕЕ
ЙПЮРЙН~$M_p^\alpha$.
мНПЛЮ Б ЩРНЛ ОПНЯРПЮМЯРБЕ НОПЕДЕКЪЕРЯЪ ЙЮЙ
$$
\|\mu\|_{M_p^\alpha}=\inf C,
$$
ЦДЕ МХФМЪЪ ЦПЮМЭ АЕПЕРЯЪ ОН БЯЕЛ ЙНМЯРЮМРЮЛ~$C$,
ДКЪ ЙНРНПШУ БШОНКМЕМН
МЕПЮБЕМЯРБН~\eqref{II:eq:negmult}.

нРЛЕРХЛ ЯКЕДСЧЫХИ ТЮЙР.
\begin{proposition}\label{II:pr:symm}
оПНЯРПЮМЯРБЮ $M^\al_p$ Х $M^\al_{p'}$ ХГНЛНПТМШ.
\end{proposition}
\begin{proof}
рЮЙ ЙЮЙ МНПЛС Б ОПНЯРПЮМЯРБЕ $H_p^{-\al}$ ЛНФМН НОПЕДЕКХРЭ ПЮБЕМЯРБНЛ
$$
 \|f\|_{-\al,p}=\supl_{\psi\in H^\al_{p'}}
\dfrac{|(f,\psi)|}{\|\psi\|_{\al,p'}},
$$
РН НЖЕМЙС~\eqref{II:eq:negmult} ЛНФМН ОЕПЕОХЯЮРЭ Б БХДЕ
$$
|(\mu(x)\phi(x),\psi(x))|\le C\|\phi\|_{\al,p}\|\psi\|_{\al,p'}
\quad \forall\phi,\psi\in D
$$
(ГДЕЯЭ ЛШ РЮЙФЕ СВХРШБЮЕЛ, ВРН $D$ ОКНРМН Б $H^\al_{p'}$).
б ОНЯКЕДМЕЕ МЕПЮБЕМЯРБН ТСМЙЖХХ $\psi$ Х $\psi$ БУНДЪР ЯХЛЛЕРПХВМН,
НРЙСДЮ ЯКЕДСЕР МЮЬЕ СРБЕПФДЕМХЕ.
\end{proof}

рЕОЕПЭ ГЮДЮДХЛЯЪ БНОПНЯНЛ: ДКЪ ЙЮЙХУ ХМДЕЙЯНБ
$\gamma\le\alpha$, $r\ge1$
ХЛЕЕР ЛЕЯРН МЕОПЕПШБМНЕ БКНФЕМХЕ Я НЖЕМЙНИ МНПЛШ
\begin{equation}\label{II:eq:incl}
H_r^{-\gamma}\subset M[H_p^\alpha\to H_p^{-\alpha}],
\qquad
\|\mu\|_{M_p^\alpha}\le C\|\mu\|_{-\gamma,r}?
\end{equation}
нРБЕР МЮ ЩРНР БНОПНЯ ДЮЕР ЯКЕДСЧЫЮЪ РЕНПЕЛЮ.

\begin{theorem}\label{II:th:negpos}
бКНФЕМХЕ~\eqref{II:eq:incl} Я НЖЕМЙНИ МНПЛ ХЛЕЕР ЛЕЯРН РНЦДЮ Х
РНКЭЙН РНЦДЮ, ЙНЦДЮ БШОНКМЪЕРЯЪ ЛСКЭРХОКХЙЮРХБМЮЪ НЖЕМЙЮ
\begin{equation}\label{II:eq:posmult}
\|\phi(x)\psi(x)\|_{\gamma,r'}
\le C\|\phi\|_{\alpha,p}\|\psi\|_{\alpha,p'}
\qquad \forall\phi,\psi\in D,
\end{equation}
ЦДЕ ВХЯКН $1\le r'\le\infty$
НОПЕДЕКЪЕРЯЪ ХГ СЯКНБХЪ ЯНОПЪФЕМХЪ
$1/r+1/r'=1$.
\end{theorem}

\begin{proof} (яП.~\cite{BSh}).
{\sl дНЯРЮРНВМНЯРЭ}.
оСЯРЭ $\mu(x)\in H_r^{-\gamma}$
Х БШОНКМЕМЮ НЖЕМЙЮ~\eqref{II:eq:posmult}. рНЦДЮ ДКЪ БЯЕУ
$\phi,\psi\in D$ ХЛЕЕЛ
\begin{align*}
\|\mu\phi\|_{-\alpha,p}
&=\sup_{0\ne\psi\in H_{p'}^\alpha}
\frac{|(\mu\phi,\psi)|}{\|\psi\|_{\alpha,p'}}
=\sup_{0\ne\psi\in H_{p'}^\alpha}
\frac{|(\mu,\ov\phi\psi)|}{\|\psi\|_{\alpha,p'}}
\\
&\le\sup_{0\ne\psi\in H_{p'}^\alpha}
\frac{\|\mu\|_{-\gamma,r}\|\ov\phi\psi\|_{\gamma,r'}}
{\|\psi\|_{\alpha,p'}}
\le C\|\mu\|_{-\gamma,r}\|\phi\|_{\alpha,p},
\end{align*}
ЦДЕ ЙНМЯРЮМРЮ~$C$ РЮ ФЕ, ВРН Х Б~\eqref{II:eq:posmult}.
оН НОПЕДЕКЕМХЧ ОНКСВЮЕЛ $\mu\in M_p^\alpha$ Х
$\|\mu\|_{M_p^\alpha}\le C\|\mu\|_{-\gamma,r}$.

{\sl мЕНАУНДХЛНЯРЭ}.
оСЯРЭ ЯОПЮБЕДКХБН БКНФЕМХЕ~\eqref{II:eq:incl}
Я НЖЕМЙНИ МНПЛ. рНЦДЮ ДКЪ $\phi,\psi\in D$
БШОНКМЕМШ МЕПЮБЕМЯРБЮ
\begin{align*}
\|\phi\psi\|_{\gamma,r'}
&=\sup_{f\in H_r^{-\gamma}}
\frac{|(f,\phi\psi)|}{\|f\|_{-\gamma,r}}
\le\sup_{\mu\in M_p^\alpha}
\frac{|(\mu\ov\phi,\psi)|}{C^{-1}\|\mu\|_{M_p^\alpha}}
\\
&\le\sup_{\mu\in M_p^\alpha}
\frac{\|\mu\ov\phi\|_{-\alpha,p}\|\psi\|_{\alpha,p'}}
{C^{-1}\|\mu\|_{M_p^\alpha}}
\le C\|\phi\|_{\alpha,p}\|\psi\|_{\alpha,p'},
\end{align*}
ОПХВЕЛ ДКЪ ОНКСВЕМХЪ ОНЯКЕДМЕЦН МЕПЮБЕМЯРБЮ ЛШ
ХЯОНКЭГНБЮКХ НЖЕМЙС~\eqref{II:eq:negmult} Я
$C=\|\mu\|_{M_p^\alpha}$.
яКЕДНБЮРЕКЭМН, БШОНКМЕМЮ НЖЕМЙЮ~\eqref{II:eq:posmult},
Х ДНЙЮГЮРЕКЭЯРБН ГЮБЕПЬЕМН.
\end{proof}

бШАЕПЕЛ МЕНРПХЖЮРЕКЭМСЧ ТСМЙЖХЧ $\eta(x)\in D$ РЮЙСЧ, ВРН
$\eta(x)\equiv1$ ОПХ $|x|\le1$. бБЕДЕЛ ОПНЯРПЮМЯРБН
$H_{p,\unif}^\alpha$
(ЯЛ.~\cite[О.~1.3.1]{MSh}),
ЯНЯРНЪЫЕЕ ХГ ТСМЙЖХИ
$u(x)\in H_{p,\loc}^\alpha$ РЮЙХУ, ВРН
$$
\|u\|_{\alpha,p,\unif}
:=\sup_{z\in\RR^n}\|\eta(x-z)u(x)\|_{\alpha,p}
<\infty.
$$
мНПЛЮ Б ОПНЯРПЮМЯРБЕ $H_{p,\unif}^\alpha$ ГЮБХЯХР НР
БШАНПЮ МЕНРПХЖЮРЕКЭМНИ ОПНАМНИ ТСМЙЖХХ $\eta(x)$,
МН ДКЪ БЯЕУ РЮЙХУ ТСМЙЖХИ ЯННРБЕРЯРБСЧЫХЕ МНПЛШ ЩЙБХБЮКЕМРМШ.
хГ НОПЕДЕКЕМХЪ КЕЦЙН ЯКЕДСЕР, ВРН
ХЛЕЧР ЛЕЯРН МЕОПЕПШБМШЕ БКНФЕМХЪ
$H^\al_p\subset H^\al_{p,\unif}$ Х
$H^\al_{p,\unif}\subset H^\al_{q,\unif}$ ОПХ $p>q$.

еЯКХ $\mu(x)\in M_p^\alpha$, РН Х
$\eta(x)\mu(x)\in M_p^\alpha$.
бБЕДЕЛ МНБСЧ МНПЛС
Б~$M_p^\alpha$:
$$
\norm\mu\norm_{M_p^\alpha}
=\sup_{z\in\RR^n}\|\eta(x-z)\mu(x)\|_{M_p^\alpha}.
$$
оПХБЕДЕММЮЪ МХФЕ РЕНПЕЛЮ ХГБЕЯРМЮ ДКЪ ЛСКЭРХОКХЙЮРНПНБ Б
ЯНАНКЕБЯЙХУ ОПНЯРПЮМЯРБЮУ Я МЕНРПХЖЮРЕКЭМШЛХ ХМДЕЙЯЮЛХ
(ЯЛ.~\cite[О.~2.1.3]{MSh}).
гДЕЯЭ ДКЪ ДНЙЮГЮРЕКЭЯРБЮ ЛШ ХЯОНКЭГСЕЛ ХДЕЧ ХГ ЯРЮРЭХ
дФ.-ц.~аЮЙЮ Х ю.~ю.~ьЙЮКХЙНБЮ~\cite{BSh}.

\begin{theorem}\label{II:th:udecomp}
мНПЛШ
$\norm\,\cdot\,\norm$ Х $\|\,\cdot\,\|$ Б~$M_p^\alpha$
ЩЙБХБЮКЕМРМШ.
\end{theorem}

\begin{proof} (яЛ.~\cite{BSh}).
мЕПЮБЕМЯРБН
$$
\norm\mu\norm_{M_p^\alpha}\le C\|\mu\|_{M_p^\alpha}
$$
ЯКЕДСЕР ХГ НОПЕДЕКЕМХЪ. дКЪ ДНЙЮГЮРЕКЭЯРБЮ НАПЮРМНИ НЖЕМЙХ
МЮЛ ОНРПЕАСЕРЯЪ ХГБЕЯРМЮЪ КЕЛЛЮ Н ПЮГАХЕМХХ ЕДХМХЖШ (ЯЛ.,
МЮОПХЛЕП,~\cite[\S 1.4]{H}).
яНЦКЮЯМН ЩРНИ КЕЛЛЕ, ЯСЫЕЯРБСЕР МЮАНП ТСМЙЖХИ
$\{\chi_j(x)\}_{j=1}^\infty$,
СДНБКЕРБНПЪЧЫХУ ЯКЕДСЧЫХЛ СЯКНБХЪЛ:
\begin{enumerate}
\item $\chi_j(x)\in D$, $0\le\chi_j(x)\le1$ Х
$\sum_{j=1}^\infty\chi_j(x)\equiv1$;
\item МНЯХРЕКХ ТСМЙЖХИ $\chi_j$ ЯНДЕПФЮРЯЪ Б ЬЮПЮУ
ПЮДХСЯЮ~$1$ (Я ПЮГКХВМШЛХ ЖЕМРПЮЛХ); ОПХ ЩРНЛ КЧАЮЪ
ТХЙЯХПНБЮММЮЪ РНВЙЮ $x\in\RR^n$ ОНЙПШБЮЕРЯЪ МНЯХРЕКЪЛХ МЕ
АНКЕЕ ВЕЛ $N$ ТСМЙЖХИ (Р.Е $\chi_j(x)\ne0$ МЕ АНКЕЕ ВЕЛ
ДКЪ $N$ ТСМЙЖХИ Х $N$ МЕ ГЮБХЯХР НР~$x$);
\item $|D^\ell\chi_j(x)|\le C$,
ЦДЕ $C$ ГЮБХЯХР КХЬЭ НР ЛСКЭРХХМДЕЙЯЮ
$\ell=(\ell_1,\dots,\ell_n)$.
\end{enumerate}
оСЯРЭ ЯЕЛЕИЯРБН $\{\chi_j\}_1^\infty$
НАКЮДЮЕР ЩРХЛХ ЯБНИЯРБЮЛХ. рНЦДЮ
$$
\sum_{j=1}^\infty\|\phi\chi_j\|_{\alpha,p}
\le C\|\phi\|_{\alpha,p},
\qquad
\phi\in D,
$$
ЦДЕ $C$ ГЮБХЯХР РНКЭЙН НР~$\alpha,p$ Х НР БШАНПЮ
$\{\chi_j\}_1^\infty$. оНЯЙНКЭЙС ОПХ БЯЕУ~$j$ МНЯХРЕКЭ
$\chi_j$ ЯНДЕПФХРЯЪ Б ЬЮПЕ ПЮДХСЯЮ~$1$, Ю $\eta(x)=1$ ОПХ
$|x|\le1$, РН ЯСЫЕЯРБСЧР $z_j\in\RR^n$ РЮЙХЕ, ВРН
$\chi_j(x)=\eta(x-z_j)\chi_j(x)$.
нАНГМЮВЮЪ $\phi_j(x)=\chi_j(x)\phi(x)$ ДКЪ $\phi\in D$,
ОНКСВЮЕЛ
\begin{multline*}
\|\mu\phi\|_{-\alpha,p'}
=\biggl\|\sum_{j=1}^\infty
\eta(x-z_j)\mu(x)\phi_j(x)\biggr\|_{-\alpha,p'}\le\\
\le\norm\mu\norm_{M_p^\alpha}
\sum_{j=1}^\infty\|\phi_j\|_{\alpha,p}
\le C\norm\mu\norm_{M_p^\alpha}\|\phi\|_{\alpha,p}.
\end{multline*}

рЕОЕПЭ ХГ НОПЕДЕКЕМХЪ ЯКЕДСЕР, ВРН
$$
\|\mu\|_{M_p^\alpha}
\le C\norm\mu\norm_{M_p^\alpha}.
$$
рЕНПЕЛЮ ДНЙЮГЮМЮ.
\end{proof}

бБЕДЕЛ НАНГМЮВЕМХЕ $\hat p=\max\{p,p'\}$. рНЦДЮ $H^\al_{\hat p,\unif}=
H^\al_{p,\unif}\cap H^\al_{p',\unif}$.

\begin{theorem}\label{II:th:descr}
бЕПМШ ЯКЕДСЧЫХЕ СРБЕПФДЕМХЪ:
\begin{enumerate}
\item  дКЪ БЯЕУ $\alpha\ge 0$  Х $p\ge 1$
ЯОПЮБЕДКХБШ БКНФЕМХЪ Я НЖЕМЙЮЛХ МНПЛ
\begin{equation}\label{II:4}
M_p^\alpha\subset H_{\hat p,\unif}^{-\alpha}, \qquad
\|\mu\|_{-\alpha,\hat p,\unif}\le C\|\mu\|_{M_p^\alpha};
\end{equation}
\item
бКНФЕМХЕ Я НЖЕМЙНИ МНПЛ
\begin{equation}\label{II:5}
H_{r,\unif}^{-\gamma}\subset M_p^\alpha, \qquad
\|\mu\|_{M_p^\alpha}\le C\|\mu\|_{-\gamma,r,\unif}
\end{equation}
ЯОПЮБЕДКХБН РНЦДЮ Х РНКЭЙН РНЦДЮ, ЙНЦДЮ БШОНКМЪЕРЯЪ
МЕПЮБЕМЯРБН~\eqref{II:eq:posmult}.
\end{enumerate}
\end{theorem}

\begin{proof}
оЕПБНЕ СРБЕПФДЕМХЕ БШРЕЙЮЕР ХГ НОПЕДЕКЕМХЪ. дЕИЯРБХРЕКЭМН,
ЕЯКХ $\mu(x)\in M_p^\alpha$ Х $\eta(x)$~-- ЦКЮДЙЮЪ
ЯПЕГЮЧЫЮЪ ТСМЙЖХЪ, РН
$$
\|\mu(x)\eta(x-z)\|_{-\alpha,p}
\le \|\mu\|_{M_p^\alpha}\|\eta(x-z)\|_{\alpha,p}
=\|\mu\|_{M_p^\alpha}\|\eta(x)\|_{\alpha,p}
=C\|\mu\|_{M_p^\alpha},
$$
ЦДЕ $C=\|\eta(x)\|_{\alpha,p}$  МЕ ГЮБХЯХР НР~$z$.
оНЩРНЛС $\mu\in H_{p,\unif}^{-\alpha}$, Х БЕПМЮ
НЖЕМЙЮ~\eqref{II:4}. оНКЭГСЪЯЭ ОПЕДКНФЕМХЕЛ~\ref{II:pr:symm}, ЮМЮКНЦХВМН
ОНКСВЮЕЛ $\mu\in H_{p',\unif}^{-\alpha}$.

бРНПНЕ СРБЕПФДЕМХЕ ЪБКЪЕРЯЪ ЯКЕДЯРБХЕЛ РЕНПЕЛ~\ref{II:th:negpos} 
Х~\ref{II:th:udecomp}.
яНЦКЮЯМН РЕНПЕЛЕ~\ref{II:th:negpos} 
Х БКНФЕМХЧ $H_r^{-\al}\subset H_{r,\unif}^{-\al}$,
НЖЕМЙЮ~\eqref{II:5} БКЕВЕР~\eqref{II:eq:posmult}. нАПЮРМН, ЕЯКХ
БШОНКМЕМЮ НЖЕМЙЮ~\eqref{II:eq:posmult} Х
$\mu\in H_{r,\unif}^{-\gamma}$, РН
$\mu(x)\eta(x-z)\in H_r^{-\gamma}\subset M_p^\alpha$ Х
$$
\|\mu(x)\eta(x-z)\|_{M_p^\alpha}
\le C\|\mu(x)\eta(x-z)\|_{-\gamma,r}.
$$
рЮЙХЛ НАПЮГНЛ,
$\norm\mu\norm_{M_p^\alpha}\le C\|\mu\|_{-\gamma,r,\unif}$,
Х~\eqref{II:5} БШРЕЙЮЕР ХГ РЕНПЕЛШ~\ref{II:th:udecomp}. рЕНПЕЛЮ ДНЙЮГЮМЮ.
\end{proof}

яНЦКЮЯМН РЕНПЕЛЕ яРПХУЮПДЯЮ (ЯЛ.~\cite{Str} ХКХ~\cite[О.~2.2.9]{MSh}),
ОПНЯРПЮМЯРБН ЛСКЭРХОКХЙЮРНПНБ
$M[H_p^\alpha\to H_p^\alpha]$ ЯНБОЮДЮЕР Я
$H_{p,\unif}^\alpha$, ЕЯКХ БШОНКЕМН СЯКНБХЕ $\alpha>n/p$.
оНКСВЕММНЕ МХФЕ СРБЕПФДЕМХЕ ЛНФМН ЯВХРЮРЭ ЮМЮКНЦНЛ ЩРНИ
РЕНПЕЛШ ДКЪ ЛСКЭРХОКХЙЮРНПНБ Б ДСЮКЭМШУ ОПНЯРПЮМЯРБЮУ.

\begin{theorem}\label{II:th:H2}
оСЯРЭ $\alpha>n/{\hat p}$. рНЦДЮ ОПНЯРПЮМЯРБЮ
$M[H_p^\alpha\to H_p^{-\alpha}]$ Х
$H_{\hat p,\unif}^{-\alpha}$ ЯНБОЮДЮЧР, Х МНПЛШ Б ЩРХУ
ОПНЯРПЮМЯРБЮУ ЩЙБХБЮКЕМРМШ.
\end{theorem}

\begin{proof}
рЮЙ ЙЮЙ ОПНЯРПЮМЯРБЮ $M^\al_p$ Х $M^\al_{p'}$
ХГНЛНПТМШ,
ДНЯРЮРНВМН ПЮЯЯЛНРПЕРЭ
ЯКСВЮИ $\hat p=p\ge2$.
бКНФЕМХЕ $M_p^\alpha\subset H_{p,\unif}^{-\alpha}$ Я НЖЕМЙНИ
МНПЛ ДКЪ БЯЕУ $\alpha\ge0$ Х $p\ge1$ ДНЙЮГЮМН Б РЕНПЕЛЕ~\ref{II:th:descr}.  
вРНАШ ДНЙЮГЮРЭ НАПЮРМНЕ БЙКЧВЕМХЕ, Б ЯННРБЕРЯРБХХ ЯН БРНПШЛ СРБЕПФДЕМХЕЛ
РЕНПЕЛШ~\ref{II:th:descr}, ДНЯРЮРНВМН ОПНБЕПХРЭ
МЕПЮБЕМЯРБН~\eqref{II:eq:posmult} ОПХ $\gamma=\alpha$ Х $r'=p$.

пЮГАЕПЕЛ НРДЕКЭМН 3 ЯКСВЮЪ:
\begin{enumerate}
\item $\al>n/{p'}$. б ЩРНЛ ЯКСВЮЕ $M[H_{p'}^\alpha\to H_{p'}^\alpha]$
ЯНБОЮДЮЕР Я $H_{p',\unif}^\alpha$ (ЯЛ.~\cite[\S\,1.7]{MSh}), ЯКЕДНБЮРЕКЭМН,
$$
\|\phi\psi\|_{\al,p'}\le \|\phi\|_{\al,p',\unif}\|\psi\|_{\al,p'}.
$$
мН $\|\phi\|_{\al,p',\unif}\le C\|\phi\|_{\al,p}$ ОПХ $p\ge p'$,
ЯКЕДНБЮРЕКЭМН, МЕПЮБЕМЯРБН~\eqref{II:eq:posmult} БШОНКМЪЕРЯЪ.

\item $\al<n/{p'}$. хЛЕЕР ЛЕЯРН МЕОПЕПШБМНЕ БКНФЕМХЕ
(ЯЛ.~\cite[\S\,2.3.1]{MSh})
$$
 H^\al_{n/\al,\unif}\cap L_\infty\subset M[H_{p'}^\alpha\to H_{p'}^\alpha].
$$
оН РЕНПЕЛЕ БКНФЕМХЪ яНАНКЕБЮ, $H^\al_p\subset L_\infty$ ОПХ $\al>n/p$, Ю
ЯОПЮБЕДКХБНЯРЭ БКНФЕМХЪ $H^\al_p\subset H^\al_{n/\al,\unif}$
ОПХ $p>n/\al$ АШКЮ НРЛЕВЕМЮ
БШЬЕ. оНЩРНЛС $\|\phi\|_{M[H_{p'}^\al\to H_{p'}^\al]} \le
\|\phi\|_\infty+\|\phi\|_{\al,n/\al,\unif}\le C\|\phi\|_{\al,p}$, Х
НЖЕМЙЮ~\eqref{II:eq:posmult} РЮЙФЕ БШОНКМЕМЮ.

\item $\al=n/{p'}$. гДЕЯЭ МЕПЮБЕМЯРБН~\eqref{II:eq:posmult} ЛНФЕР АШРЭ
ОНКСВЕМН ОПХЛЕМЕМХЕЛ ОНКХКХМЕИМНИ ХМРЕПОНКЪЖХНММНИ РЕНПЕЛШ
(ЯЛ., МЮОПХЛЕП,~\cite[\S\,4.4]{BL})
ОН ОЮПЮЛЕРПС $\al$
Й ОПЕДШДСЫХЛ ДБСЛ ЯКСВЮЪЛ.

\end{enumerate}
рЕНПЕЛЮ ДНЙЮГЮМЮ.
\end{proof}

б ЯКСВЮЕ $\alpha\le n/{\hat p}$ ДЮРЭ НОХЯЮМХЕ ОПНЯРПЮМЯРБЮ
$M_p^\alpha$ Б РЕПЛХМЮУ ОПНЯРПЮМЯРБ яНАНКЕБЮ СФЕ
МЕБНГЛНФМН. оН ЮМЮКНЦХХ Я РЕНПХЕИ ЛСКЭРХОКХЙЮРНПНБ Б
ЯНАНКЕБЯЙХУ ОПНЯРПЮМЯРБЮУ Я МЕНРПХЖЮРЕКЭМШЛХ ХМДЕЙЯЮЛХ
ЦКЮДЙНЯРХ, РЮЙНЕ НОХЯЮМХЕ ЛНФМН ОПНБЕЯРХ Б РЕПЛХМЮУ
ЕЛЙНЯРХ, УНРЪ ЩРН РПЕАСЕР ДНОНКМХРЕКЭМНИ ЯЕПЭЕГМНИ ПЮАНРШ
(ЯЛ ПЕГСКЭРЮРШ бЕПАХЖЙНЦН Б~\cite{MSh}). нДМЮЙН ОПНБЕПЙЮ
СЯКНБХИ Б РЕПЛХМЮУ ЕЛЙНЯРЕИ РПСДМН НЯСЫЕЯРБХЛЮ, ОНЩРНЛС
ОНКЕГМН ХЛЕРЭ ЩТТЕЙРХБМШЕ ДНЯРЮРНВМШЕ СЯКНБХЪ
ОПХМЮДКЕФМНЯРХ ТСМЙЖХИ ПЮЯЯЛЮРПХБЮЕЛШЛ ОПНЯРПЮМЯРБЮЛ
ЛСКЭРХОКХЙЮРНПНБ. мЮЛ СДЮКНЯЭ МЮИРХ РНВМШЕ СЯКНБХЪ МЮ
ХМДЕЙЯШ $\gamma$ Х~$r$, ОПХ ЙНРНПШУ ЯОПЮБЕДКХБН БКНФЕМХЕ
$H^{-\gamma}_{r,\unif}\subset M^\alpha_p$
БЛЕЯРЕ Я НЖЕМЙНИ МНПЛ. гДЕЯЭ ЛШ ОПХБЕДЕЛ ДБЕ РЕНПЕЛШ РЮЙНЦН
РХОЮ. пЕГСКЭРЮР ОЕПБНИ ХГ МХУ БЮФЕМ КХЬЭ ОПХ $\alpha=n/p$,
ОНЯЙНКЭЙС ОПХ $\alpha<n/p$ АНКЕЕ РНВМНЕ СРБЕПФДЕМХЕ ДЮЕР
БРНПЮЪ ХГ ЩРХУ РЕНПЕЛ.

\begin{theorem}\label{II:th:Hp}
оСЯРЭ $\al\le n/p$.
рНЦДЮ ЯОПЮБЕДКХБН БКНФЕМХЕ
$H^{-\gamma}_{r,\unif}\subset M^\alpha_p$
БЛЕЯРЕ Я НЖЕМЙНИ МНПЛ, ЕЯКХ ВХЯКЮ $\gamma$ Х~$r$
ОНДВХМЕМШ СЯКНБХЪЛ
$$
\gamma\le\alpha
\qquad \text{Х}\qquad
r>\max\{1,\frac{n}{2\al-\ga}\}.
$$
\end{theorem}

\begin{proof}
дНЯРЮРНВМН ДНЙЮГЮРЭ ЩРН СРБЕПФДЕМХЕ ОПХ $\gamma=\alpha$.
дЕИЯРБХРЕКЭМН, ОПЕДОНКНФХЛ, ВРН ДНЙЮГЮМН БКНФЕМХЕ
$H^\alpha_{r_1}\subset M^\alpha_p$
БЛЕЯРЕ Я НЖЕМЙНИ МНПЛ ОПХ
$r_1>n/\alpha$.
б ЯХКС РЕНПЕЛШ яНАНКЕБЮ,
$H^{-\gamma}_r\subset H^{-\alpha}_{r_1}$,
ЕЯКХ
$\alpha-\gamma= n/r-n/r_1\ge0$ Х $r>1$.
яКЕДНБЮРЕКЭМН,
$H^\gamma_r\subset M^\alpha_p$,
ЕЯКХ
$r= \dfrac{r_1}{(1+r_1(\alpha-\gamma)/n)}$,
Р.Е ЕЯКХ
$r> n/(2\alpha-\gamma)$.

хРЮЙ, ОСЯРЭ $\gamma=\alpha$.
оНКНФХЛ $\nu=n/p+\ep$, $\ep>0$,
Х ГЮОХЬЕЛ ЯКЕДСЧЫХЕ МЕПЮБЕМЯРБЮ:
$$
\|\phi\psi\|_{\nu,p'}\le C\|\phi\|_{\nu,p}\|\psi\|_{\nu,p'},
\qquad
\|\phi\psi\|_{0,1}\le\|\phi\|_{0,p}\|\psi\|_{0,p'}.
$$
оЕПБНЕ ХГ МХУ ДНЙЮГЮМН Б ОПЕДШДСЫЕИ РЕНПЕЛЕ (ДКЪ $\nu>n/p$),
Ю БРНПНЕ ЕЯРЭ НАШВМНЕ
МЕПЮБЕМЯРБН ц╦КЭДЕПЮ. рЕОЕПЭ, ХЯОНКЭГСЪ ОНКХКХМЕИМСЧ ХМРЕПНКЪЖХНММСЧ
РЕНПЕЛС (ЯЛ., МЮОПХЛЕП,~\cite[\S\,4.4]{BL}), ЛШ ОПХУНДХЛ Й МЕПЮБЕМЯРБС $$
\|\phi\psi\|_{\alpha,s}
\le C\|\phi\|_{\alpha,p}\|\psi\|_{\alpha,p'}
$$
ДКЪ
$$
0\le\alpha\le\nu,
\qquad
s=\frac{n+p\ep}{n+p\ep-\al},
\qquad
s'=\frac{n+p\ep}{\al}.
$$
дКЪ ГЮБЕПЬЕМХЪ ДНЙЮГЮРЕКЭЯРБЮ МЮДН БШАПЮРЭ
$\ep$ РЮЙХЛ, ВРНАШ $s'=r$, Ю ГЮРЕЛ
БНЯОНКЭГНБЮРЭЯЪ РЕНПЕЛЮЛХ~\ref{II:th:negpos} Х~\ref{II:th:descr}.
\end{proof}

б ЯКСВЮЕ $\alpha=n/p$ ПЕГСКЭРЮР МЕКЭГЪ СКСВЬХРЭ, ГЮЛЕМХБ
МЕПЮБЕМЯРБН ДКЪ~$r$ МЮ ЯННРБЕРЯРБСЧЫЕЕ ПЮБЕМЯРБН (Р.Е
ОНКНФХРЭ $\ep=0$). нДМЮЙН ЩРН ЛНФМН ЯДЕКЮРЭ ОПХ
$\alpha<n/p$. дНЙЮГЮРЕКЭЯРБН ЩРНЦН ТЮЙРЮ РПЕАСЕР АНКЕЕ
ЦКСАНЙНЦН ЮМЮКХГЮ. нРЛЕРХЛ, ВРН ОПХ ЖЕКНЛ~$\alpha$
ДНЙЮГЮРЕКЭЯРБН ЛНФМН ОПНБЕЯРХ ГМЮВХРЕКЭМН ОПНЫЕ
(ЯЛ.~\cite{NeSh}).

\begin{theorem}\label{II:th:Polking}
оСЯРЭ $p>1$ Х $\alpha<n/p$.
рНЦДЮ ЯОПЮБЕДКХБН БКНФЕМХЕ
$H^{-\gamma}_{r,\unif}\subset M^\alpha_p$
БЛЕЯРЕ Я НЖЕМЙНИ МНПЛ, ЕЯКХ ВХЯКЮ $\gamma,r$
ОНДВХМЕМШ СЯКНБХЪЛ
$$
\gamma\le\alpha
\qquad\text{Х}\qquad
r\ge\frac{n}{2\alpha-\gamma},\qquad r>1.
$$
\end{theorem}

\begin{proof}
йЮЙ Х Б ОПЕДШДСЫЕИ РЕНПЕЛЕ, ЛНФМН НЦПЮМХВХРЭЯЪ ЯКСВЮЕЛ
$\gamma=\alpha$,
ОПХВЕЛ СРБЕПФДЕМХЕ МСФМН ДНЙЮГЮРЭ РНКЭЙН ДКЪ
$r=n/{\alpha}$.
яНЦКЮЯМН РЕНПЕЛЕ~\ref{II:th:descr}, МСФМН ДНЙЮГЮРЭ НЖЕМЙС
\begin{equation}\label{II:6}
\|\phi\psi\|_{\alpha,s}
\le C\|\phi\|_{\alpha,p}\|\psi\|_{\alpha,p'}
\qquad \forall\phi,\psi\in D,
\end{equation}
ЦДЕ $s=n/(n-\alpha)>1$.
дКЪ ДНЙЮГЮРЕКЭЯРБЮ ЩРНИ НЖЕМЙХ
ЛШ АСДЕЛ ХЯОНКЭГНБЮРЭ ЛЕРНДШ, ПЮГБХРШЕ Б ЯРЮРЭЕ
оНКЙХМЦЮ~\cite{Pol}.

пЮЯЯЛНРПХЛ МЕКХМЕИМШИ НОЕПЮРНП
$$
D_{\rho,q}^\eta u(x)
=\int_0^\infty\biggl(\int_{B_1}
|u(x+\xi y)-u(x)|^\rho\,dy\biggr)^{q/\rho}\xi^{-1-\eta q}\,d\xi,
$$
ЦДЕ $\eta>0$, $\rho\ge1$, $q\ge2$.
мЮЛ ОНРПЕАСЧРЯЪ
ЯКЕДСЧЫХЕ СРБЕПФДЕМХЪ, ДНЙЮГЮММШЕ Б ПЮАНРЮУ
оНКЙХМЦЮ~\cite{Pol} Х яРПХУЮПДЯЮ~\cite{Str}.

\begin{theorem*}{оПЕДКНФЕМХЕ A \rm(ЯЛ.~\cite{Pol})}
оСЯРЭ $0<\gamma<1$, $1\le\rho\le q$, $2\le q<\infty$,
$r>\max\{1, n\rho/(n+\gamma\rho)\}$.
рНЦДЮ БШОНКМЕМН ЯКЕДСЧЫЕЕ МЕПЮБЕМЯРБН:
$$
\|D_{\rho,q}^\ga u(x)\|_{0,r}\le C\|u(x)\|_{\ga,r},
\qquad
u\in H_r^\eta,
$$
ЦДЕ ЙНМЯРЮМРЮ $C$ ГЮБХЯХР НР~$\ga,\rho,q,r,n$,
МН МЕ ГЮБХЯХР НР~$u(x)$.
\end{theorem*}

\begin{theorem*}{оПЕДКНФЕМХЕ B \rm(ЯЛ.~\cite{Pol})}
оСЯРЭ $1/\rho_1+1/\rho_2=1/\rho$, $1/q_1+1/q_2=1/q$ Х
$\eta_1+\eta_2=\eta$. рНЦДЮ ДКЪ БЯЕУ
$u\in L_{\rho_1,\loc}$ Х $v\in L_{\rho_2,\loc}$
БШОНКМЕМН МЕПЮБЕМЯРБН
$$
|D_{\rho,q}^\eta(uv)(x)|
\le|u(x)D_{\rho,q}^\eta v(x)|
+|v(x)D_{\rho,q}^\eta u(x)|
+|D_{\rho_1,q_1}^{\eta_1}u(x)D_{\rho_2,q_2}^{\eta_2}v(x)|.
$$
\end{theorem*}

\begin{theorem*}{оПЕДКНФЕМХЕ C \rm(ЯЛ.~\cite{Pol,Str})}
хЛЕЕР ЛЕЯРН ЯКЕДСЧЫЕЕ ЯННРМНЬЕМХЕ:
$$
\|v\|_{\eta,p}\asymp\|D_{1,2}^\eta v\|_{0,p}+\|v\|_{0,p},
\qquad p\ge 1,\ \eta >0, \qquad v(x)\in H^\eta_p,
$$
ЦДЕ ЯХЛБНК $\asymp$ НАНГМЮВЮЕР ЩЙБХБЮКЕМРМНЯРЭ (Р.Е НРМНЬЕМХЕ
ОПЮБНИ Х КЕБНИ ВЮЯРЕИ НЦПЮМХВЕМН Я НАЕХУ ЯРНПНМ
ЙНМЯРЮМРЮЛХ, МЕ ГЮБХЯЪЫХЛХ НР $v(x)$).
\end{theorem*}

хЯОНКЭГСЪ ЩРХ СРБЕПФДЕМХЪ, ДНЙЮФЕЛ НЖЕМЙС~\eqref{II:6} ДКЪ
$s=n/(n-\alpha)$. оНКНФХЛ
$\alpha=[\alpha]+\eta$, ЦДЕ $[\alpha]$
НАНГМЮВЮЕР ЖЕКСЧ
ВЮЯРЭ~$\alpha$ Х $0\le\eta<1$. аСДЕЛ ОПЕДОНКЮЦЮРЭ, ВРН
 $\eta>0$ (ДКЪ
$\eta=0$ ДНЙЮГЮРЕКЭЯРБН  ОПНЫЕ Х ЛНФЕР АШРЭ
ОНКСВЕМН ОПХ ОНЛНЫХ НАШВМНИ ТНПЛСКШ кЕИАМХЖЮ).
уНПНЬН ХГБЕЯРМН ЯКЕДСЧЫЕЕ ЯННРМНЬЕМХЕ
(ЯЛ., МЮОПХЛЕП,~\cite{Ste}):
$$
\|\phi\psi\|_{\alpha,s}
\asymp\sum_{|\ell|\le[\alpha]}\|D^\ell(\phi\psi)\|_{\eta,s}.
$$
яНЦКЮЯМН ОПЕДКНФЕМХЧ~C, ОПЮБЮЪ ВЮЯРЭ НЖЕМХБЮЕРЯЪ ЯБЕПУС
БЕКХВХМНИ (Б ДЮКЭМЕИЬЕЛ Б НЖЕМЙЮУ ЛШ НОСЯЙЮЕЛ ЙНМЯРЮМРШ)
$$
\sum_{|j|+|m|\le[\alpha]}
\bigl(\|D_{1,2}^\eta(D^j\phi D^m\psi)\|_{0,s}
+\|D^j\phi D^m\psi\|_{0,s}\bigr).
$$
щРН БШПЮФЕМХЕ ЛНФЕР АШРЭ НЖЕМЕМН ЯБЕПУС, Я СВЕРНЛ
ОПЕДКНФЕМХЪ~B, ЯКЕДСЧЫЕИ БЕКХВХМНИ (ЛШ АСДЕЛ НЖЕМХБЮРЭ
РНКЭЙН ОЕПБНЕ ЯКЮЦЮЕЛНЕ, БРНПНЕ НЖЕМХБЮЕРЯЪ ЯСЫЕЯРБЕММН ОПНЫЕ):
\begin{align}
\sum_{|j|+|m|\le[\alpha]}
&
\bigl(\|D^m\psi D_{1,2}^\eta(D^j\phi)\|_{0,s}
+\|D^j\phi D_{1,2}^\eta(D^m\psi)\|_{0,s}\notag
\\ &\qquad
+\|(D_{\rho_1,q_1}^{\eta_1}D^j\phi)
(D_{\rho_2,q_2}^{\eta_2}D^m\psi)\|_{0,s}\bigr).
\label{II:7}
\end{align}
рЕОЕПЭ БШАЕПЕЛ ОНДУНДЪЫХЛ НАПЮГНЛ
ВХЯКЮ~$\eta_{1,2}$, $\rho_{1,2}$, $q_{1,2}$. оНКНФХЛ
\begin{align}
&0<\eta_1<\max\biggl\{\eta,\frac np-\alpha\biggr\},
\qquad
\eta_2=\eta-\eta_1,
\label{II:8}
\\
&\rho_1=\frac {p(n-\al)}{n+p(|j|+\eta_1-\al)},
\quad
\rho_2=\frac {p(n-\al)}{pn-n-p(|j|+\eta_1)},
\label{II:9}
\\
&q_1=2\rho_1,
\quad
q_2=2\rho_2.
\end{align}

рЮЙ ЙЮЙ $|j|+\eta_1<\alpha$, РН $\rho_1,\rho_2>1$ Х
$1/\rho_1+1/\rho_2=1$.

нЖЕМХЛ РПЕРЭЕ ЯКЮЦЮЕЛНЕ Б~\eqref{II:7} Я 
ОНЛНЫЭЧ МЕПЮБЕМЯРБЮ ц╦КЭДЕПЮ БЕКХВХМНИ
$$
\bigl(\|D_{\rho_1,2\rho_1}^{\eta_1}D^j\phi\|_{0,\rho_1s}\bigr)
\bigl(\|D_{\rho_2,2\rho_2}^{\eta_2}D^m\psi\|_{0,\rho_2s}\bigr).
$$
яНЦКЮЯМН ОПЕДКНФЕМХЧ~A,
ОНЯКЕДМЕЕ БШПЮФЕМХЕ ЛНФМН НЖЕМХРЭ БЕКХВХМНИ
\begin{equation}\label{II:10}
\|D^j\phi\|_{\eta_1,\rho_1s}\|D^m\psi\|_{\eta_2,\rho_2s},
\end{equation}
ЕЯКХ РНКЭЙН БШОНКМЕМШ СЯКНБХЪ
$$
\rho_1 s>1,
\quad
\rho_1 s>\frac{n\rho_1}{n+\eta_1\rho_1},
\qquad
\rho_2 s>1,
\quad
\rho_2 s>\frac{n\rho_2}{n+\eta_2\rho_2},
$$
ЯОПЮБЕДКХБНЯРЭ ЙНРНПШУ ЯКЕДСЕР ХГ РНЦН, ВРН $s>1$ Х $\rho_{1,2}>1$.

дЮКЕЕ, БШПЮФЕМХЕ~\eqref{II:10} НЖЕМХБЮЕРЯЪ БЕКХВХМНИ
$$
\|\phi\|_{|j|+\eta_1,\rho_1 s}\|\psi\|_{|m|+\eta_2,\rho_2 s}
\le C\|\phi\|_{\alpha,p}\|\psi\|_{\alpha,p'}.
$$
оНЯКЕДМЪЪ НЖЕМЙЮ ОНКСВЮЕРЯЪ ХГ РЕНПЕЛШ БКНФЕМХЪ яНАНКЕБЮ Х
МЕПЮБЕМЯРБ
$$
|j|+\eta_1-\frac n{\rho_1 s}\le\alpha-\frac np,
\qquad
|m|+\eta_2-\frac n{\rho_2 s}\le\alpha-\frac n{p'}.
$$
оЕПБНЕ ХГ ЩРХУ МЕПЮБЕМЯРБ ЯКЕДСЕР ХГ~\eqref{II:9} Х СЯКНБХЪ
$|j|+\eta_1<\alpha$, Ю БРНПНЕ~-- ХГ~\eqref{II:9} Х СЯКНБХЪ
$|m|+\eta_2<\alpha$.

нЖЕМЙХ ОЕПБНЦН Х БРНПНЦН ЯКЮЦЮЕЛНЦН Б~\eqref{II:7}
ЛНЦСР АШРЭ ОНКСВЕМШ РЮЙХЛ ФЕ НАПЮГНЛ, ЕЯКХ ОНКНФХРЭ
$\eta_1=\eta$, $\eta_2=0$ ДКЪ ОЕПБНЦН ЯКЮЦЮЕЛНЦН Х
$\eta_1=0$, $\eta_2=\eta$ ДКЪ БРНПНЦН.
рЕНПЕЛЮ ДНЙЮГЮМЮ.
\end{proof}

йЮЙ ЯКЕДЯРБХЕ ОНКСВЮЕЛ ЯКЕДСЧЫЕЕ СРБЕПФДЕМХЕ.
\begin{theorem}\label{II:th:main}
оСЯРЭ $m>n/{\hat p}$. рНЦДЮ ДКЪ РНЦН, ВРНАШ НОЕПЮРНП
$$
  Lu(x)=[(-\De)^m+q(x)]u(x)
$$
АШК НЦПЮМХВЕМ ЙЮЙ НОЕПЮРНП ХГ $H^m_p$ Б $H^{-m}_p$, МЕНАУНДХЛН
Х ДНЯРЮРНВМН, ВРНАШ
$$
  q\in H^{-m}_{\hat p,\unif}.
$$
оПХ $m\le n/{\hat p}$ ДНЯРЮРНВМШЛ СЯКНБХЕЛ НЦПЮМХВЕММНЯРХ НОЕПЮРНПЮ $L$
ХГ $H^m_p$ Б $H^{-m}_p$ ЪБКЪЕРЯЪ СЯКНБХЕ
$$
  q\in H^{-m}_{r,\unif},\quad r>n/m.
$$
еЯКХ $m<n/{\hat p}$ Х $p>1$, РН ЛНФМН ОНКНФХРЭ $r=n/m$.
\end{theorem}

\begin{proof}
нОЕПЮРНП $(-\De)^m$ НЦПЮМХВЕМ ХГ $H^m_p$ Б $H^{-m}_p$. нЦПЮМХВЕММНЯРЭ НОЕПЮРНПЮ
СЛМНФЕМХЪ МЮ $q$ ХГ $H^m_p$ Б $H^{-m}_p$ ЯКЕДСЕР ХГ РЕНПЕЛ~\ref{II:th:H2}-
\ref{II:th:Polking}.
\end{proof}

б ЯКСВЮЕ НДМНЛЕПМНЦН ОПНЯРПЮМЯРБЮ ($n=1$), Ю РЮЙФЕ $m=1$, ЮМЮКНЦХВМШИ ПЕГСКЭРЮР НА 
НЦПЮМХВЕММНЯРХ НОЕПЮРНПЮ $L$, МН Б ДПСЦНИ ТНПЛЕ, ОНКСВЕМ Б МЕДЮБМЕИ ПЮАНРЕ
б.~ц.~лЮГЭХ Х х.~е.~бЕПАХЖЙНЦН~\cite{MV}.

нЯМНБМШЕ ПЕГСКЭРЮРШ ЩРНИ ПЮАНРШ АШКХ ХГКНФЕМШ Б ДХЯЯЕПРЮЖХХ 
л.~х.~мЕИЛЮМ--ГЮДЕ~\cite{Nth}.

\end{document}